\documentclass{elsart}
\usepackage{ifpdf}
\usepackage{graphicx,amssymb,lineno}
\ifpdf
\usepackage[%
  pdftitle={Instructions for use of the document class
    elsart},%
  pdfauthor={Simon Pepping},%
  pdfsubject={The preprint document class elsart},%
  pdfkeywords={instructions for use, elsart, document class},%
  pdfstartview=FitH,%
  bookmarks=true,%
  bookmarksopen=true,%
  breaklinks=true,%
  colorlinks=true,%
  linkcolor=blue,anchorcolor=blue,%
  citecolor=blue,filecolor=blue,%
  menucolor=blue,pagecolor=blue,%
  urlcolor=blue]{hyperref}
\else
\usepackage[%
  breaklinks=true,%
  colorlinks=true,%
  linkcolor=blue,anchorcolor=blue,%
  citecolor=blue,filecolor=blue,%
  menucolor=blue,pagecolor=blue,%
  urlcolor=blue]{hyperref}
\fi

\makeatletter
\def\elsartstyle{%
    \def\normalsize{\@setfontsize\normalsize\@xiipt{14.5}}
    \def\small{\@setfontsize\small\@xipt{13.6}}
    \let\footnotesize=\small
    \def\large{\@setfontsize\large\@xivpt{18}}
    \def\Large{\@setfontsize\Large\@xviipt{22}}
    \skip\@mpfootins = 18\p@ \@plus 2\p@
    \normalsize
} \@ifundefined{square}{}{} \makeatother

\usepackage{amstext}
\usepackage{amssymb}
\usepackage{graphicx}
\usepackage{amsmath}
\usepackage{epsf}
\usepackage{graphpap}
\def\be{\begin{equation}}
\def\ee{\end{equation}}
\begin{document}
\begin{frontmatter}

\title{The applications of the partial Hamiltonian approach to mechanics and other areas}
\author{R. Naz$^*$}
\address{Centre for Mathematics and Statistical Sciences,
  Lahore School of Economics, Lahore, 53200, Pakistan\\
  $^*$  Corresponding Author Email: drrehana@lahoreschool.edu.pk.\\
   Tel:
 92-3315439545\\
}
  \begin{abstract}
  The partial Hamiltonian systems of the form $\dot q^i=\frac{\partial
H}{\partial p_i},
 \dot p^i=-\frac{\partial H}{\partial q_i}+\Gamma^i(t,q^i,p_i)$ arise widely in different fields of the applied mathematics. The
partial Hamiltonian systems appear for a mechanical system with
non-holonomic nonlinear constraints and non-potential generalized
forces. In dynamic optimization problems of economic growth theory
involving a non-zero discount factor the partial Hamiltonian systems
arise and are known as a current value Hamiltonian systems.
  It is shown that the partial Hamiltonian approach proposed
earlier for the current value Hamiltonian systems arising in
economic growth theory  Naz et al \cite{naz}  is applicable to
mechanics and other areas as well. The partial Hamiltonian approach
is utilized to construct first integrals and closed form solutions
of optimal growth model with environmental asset, equations of
motion for a mechanical system with non-potential forces, the
force-free Duffing Van Der Pol Oscillator and Lotka-Volterra models.
\end{abstract}


\end{frontmatter}
\section{Introduction}
The classical mechanics was reformulated as Lagrangian mechanics by
Joseph Louis Lagrange in 1788 \cite{lag}. In 1833, William Rowan
Hamilton \cite{ham} formulated the Hamiltonian mechanics by
utilizing the Legendre transformation \cite{leg,leg2}. Later on, the
notions of Lagrangian and Hamiltonian became popular in other fields
as well e.g. continuum mechanics, fluid mechanics,  quantum
mechanics, plasma physics, engineering, mathematical biology,
economic growth theory and many other fields dealing with dynamic
optimization problems. A dynamic optimization problem involves the
determination of the extremal of the functional involving time,
dependent, independent variables and their derivative up to finite
order. There are three major approaches to deal with dynamic
optimization problems: calculus of variations, dynamic programming
and optimal control theory. The calculus of variation utilizes the
notion of a standard Lagrangian and provides a set of equations
known as Euler-Lagrange equations \cite{lag,eul}. The dynamic
programming was introduced by Richard Ernest Bellman \cite{bel}. The
optimal control theory is an extension of calculus of variation and
is developed by Lev Semyonovich Pontryagin \cite{pont}.

 The first integrals or conservation laws for differential equations are
essential in constructing exact solutions (see
e.g.\cite{nlm1}-\cite{nlm3} and references therein).  The first
integrals/conservation laws for the Euler-Lagrange differential
equations can be established with the help of celebrated Noether's
theorem \cite{Noe} provided a standard Lagrangian exists. Most of
the differential equations that describe the real world phenomena do
not admit standard Lagrangian and thus Noether's theorem \cite{Noe}
cannot be applied to construct the first integrals and conservation
laws. The partial Lagrangian approach \cite{imran,kar3} was
developed to construct first integrals/conservation laws for
differential equations which do not have standard Lagrangian. The
partial or discount free Lagrangian approach \cite{cl} is developed
to derive the first integrals and closed-form solutions for the
calculus of variation problems involving partial or discount free
Lagrangian in economic growth theory. Naeem and Mahomed \cite{nae}
provided notions of approximate partial Lagrangian and approximate
Euler-Lagrange equations for perturbed ODEs. There are methods to
obtain first integrals which do not rely on the knowledge of a
Lagrangian function. The characteristic method \cite{ste} and direct
method \cite{olv,anc} have been successfully applied to establish
first integrals of several differential equations. The most
effective and systematic Maple based computer package GeM developed
by Cheviakov  \cite{gem1}-\cite{gem2} works in an excellent way to
derive first integrals. A review of all different approaches to
construct first integrals for differential equations is presented in
\cite{comp1,comp2}.

  The
Legendre transformation provides the equivalence of the
Euler-Lagrange and Hamiltonian equations \cite{leg2}. In 2010,
Dorodnitsyn and Kozlov \cite{rom} established the relation between
symmetries and first integrals for both continuous canonical
Hamiltonian equations and discrete Hamiltonian equations by
utilizing the Legendre transformation. Thus well-known Noether's
theorem was formulated in terms of the Hamiltonian function and
symmetry operators. A current value Hamiltonian approach was
proposed by Naz et al \cite{naz,naz2016} to derive the first
integrals and closed-form solutions for systems of first-order ODEs
arising from the optimal control problems involving current value
Hamiltonian in economic growth theory. Mahomed and Roberts
\cite{diff} focused on the characterization of Hamiltonian
symmetries and their first integrals. This is applicable to standard
Hamiltonian systems.

  In this paper, I focus on the partial Hamiltonian systems of the form $\dot q^i=\frac{\partial
H}{\partial p_i},
 \dot p^i=-\frac{\partial H}{\partial q_i}+\Gamma^i(t,q^i,p_i)$ which arises widely in economic growth
theory, physics, mechanics, biology and in some other fields of
applied mathematics. The partial Hamiltonian systems appear for the
mechanical systems with non-holonomic nonlinear constraints and
non-potential generalized forces. In dynamic optimization problems
of economic growth theory involving a non-zero discount factor the
partial Hamiltonian systems arise and are known as  current value
Hamiltonian systems. These type of systems arise also in classical
field theories almost every conservative fluid and plasma theory has
this form. The partial Hamiltonian systems have a real physical
structure.

The layout of the paper is as follows. Preliminaries on known forms
of partial Hamiltonian systems are presented in Section 2. In
Section 3, real world applications are presented to show
effectiveness of partial Hamiltonian approach. The first integrals
of the optimal growth model with environmental asset, equations of
motion for a mechanical system, the force-free duffing Van Der Pol
oscillator and Lotka-Volterra system are established. Finally,
conclusions are presented in Section 4.

\section{Preliminaries}

 Let $t$ be the independent variable which
is usually time and $(q,p)=(q^1,...,q^n, p_1,...,p_n)$ the phase
space coordinates.  The following results are adopted from
\cite{rom, naz, naz2016}:
 The Euler operator ${\delta}/{\delta
q^i}$ and the variational operator ${\delta}/{\delta p_i}$ are
defined as \be \frac{\delta}{\delta q^i}=\frac{\partial}{\partial
q^i}-D\frac{\partial}{\partial \dot q^i},\; i=1,2, \cdots, n,
\label{(1)} \ee and \be \frac{\delta}{\delta
p_i}=\frac{\partial}{\partial p_i}-D\frac{\partial}{\partial \dot
p_i}, i=1,2,\; \cdots, n, \label{(2)}\ee where \be
D_t=\frac{\partial }{\partial t}+\dot q^i\frac{\partial }{\partial
q^i}+\dot p_i\frac{\partial }{\partial p_i}+\cdots\label{(2t)} \ee
is the total derivative operator with respect to the time $t$. The
summation convention applies for repeated indices here and in the
sequel.

The variables $t, q^i, p_i$ are independent and connected by the
differential relations \be \dot p_i=D_t(p_i),\; \dot
q^i=D_t(q^i),\;i=1,2,\cdots,n. \label{(3)} \ee

A partial Hamiltonian $H$ satisfies (see \cite{naz})
\begin{eqnarray}
\dot q^i=\frac{\partial H}{\partial p_i},
\nonumber\\
[-1.5ex]\label{(6rn)} \\[-1.5ex]
 \dot p^i=-\frac{\partial H}{\partial q^i}+\Gamma^i(t,q^i,p_i),\; i=1,\ldots,n,  \nonumber
\end{eqnarray}
 where $\Gamma^i$ are in general non-zero functions of $t,\;q^i,\;p^i$.

  The generators of
point symmetries in the space $t,q ,p$ are operators of the form
\cite{rom,naz,naz2016}
 \be X=\xi(t,q,p)\frac{\partial }{\partial t}+ \eta ^i (t, q,p)\frac{\partial }{\partial q^i}
 + \zeta _i (t,q,p) \frac{\partial }{\partial p_i}. \label{(7)}
 \ee

Naz et al \cite{naz,naz2016} provided following criteria to derive
partial Hamiltonian operators and associated first integrals:

  An operator $X$ of the form (\ref{(7)})
is said to be a partial Hamiltonian operator corresponding to a
current value Hamiltonian $H(t,q,p)$ , if there exists a function
$B(t,q,p)$ such that
 \be
 \zeta_i \frac{\partial H}{\partial p_i}+p_i D_t(\eta^i)-X(H)
-H D_t(\xi)=D_t(B) +(\eta^i -\xi \frac{\partial H}{\partial
p_i})(-\Gamma^i) \label{(9)} \ee holds on the system (\ref{(6rn)})
then system has a first-integral \be I= p_i \eta ^i -\xi H-B.
\label{(10)} \ee

The Hamiltonian system of form (\ref{(6rn)}) arises in economic
growth theory, mechanics, physics, biology and in some other fields
of applied mathematics.

{\bf Remark 1:} The function $H$ which results in a Hamiltonian
system of form (\ref{(6rn)}) is defined in different ways. In
Economic growth theory $H$ is a current value Hamiltonian function
see e.g. \cite{naz, naz2016, chiang}. In a first-order mechanical
system with non-holonomic constraints or non-potential generalized
forces or external forces the function $H$ which gives rise to a
standard Hamiltonian but it yields a Hamiltonian system of form
(\ref{(6rn)}) see e.g. \cite{mec,mec1}. One can also find these
structures in different fields of applied Mathematics. These systems
have a real physical structure.

A natural question arises what is significance of functions
$\Gamma^i(t,q^i,p_i)$?

{\bf Remark 2:} In each field the functions $\Gamma^i(t,q^i,p_i)$
are interpreted in different ways but are always some physical
quantities. In correspondence to the economic growth theory
$\Gamma^i$ is associated with discount factor. In mechanics the
systems of form  {(\ref{(6rn)})} arise and the functions $\Gamma^i$
contain non-potential generalized forces. In some other mechanical
systems the functions $\Gamma^i$ are related to the generalized
constrained forces. In other fields as well the functions
$\Gamma^i(t,q^i,p_i)$ have structural properties and describe some
physical phenomena.

\section{Applications}
In this Section, the first-integrals and closed-form solutions for
some real world models are derived to show effectiveness of partial
Hamiltonian approach developed by Naz et al \cite{naz}.  The first
integrals and closed-form solutions of optimal growth model with
environmental asset are established.  The first integrals of the
equations of motion for a mechanical system are established. Both of
these models have real Hamiltonian structure. The force-free Duffing
Van Der Pol Oscillator and Lotka-Volterra model have no real
Hamiltonian structure but these can be expressed as a partial
Hamiltonian system.
\subsection{Optimal growth model with environmental asset}
 The partial Hamiltonian approach \cite{naz} is applied to
the optimal growth model with an environmental assets investigated
by Le Kama and Schubert \cite{lek}. The social planner seeks to
maximize \be Max \quad
\int_0^{\infty}\frac{(cs^{\phi})^{1-\sigma}}{1-\sigma} e^{-{\rho t}
} dt ,\; \phi \geq 0, \sigma>0, \sigma\not=1, \rho>0,\label{(e1)}\ee
subject to the constraint \be \dot s =m s -c,\label{(e2)}\ee where
$s(t)$ is the stock of environmental asset, $c(t)$ is the
consumption, $\sigma$ is the inverse of intertemporal elasticity of
substitution, $\phi$ is the relative preference for environment,
$m>0$ is the regeneration rate for environmental asset and $\rho$ is
the discount factor.
 \subsubsection{First integrals} The current value
Hamiltonian for this model is \be
H(t,c,s,p)=\frac{(cs^{\phi})^{1-\sigma}}{1-\sigma} +p(m s -c)
\label{(e3)},\ee where $p(t)$  is the costate variable. The
necessary first order conditions for optimal control yield
 \be
p=c^{-\sigma} s^{\phi(1-\sigma)}, \label{(e4)} \ee
 \be \dot s
=m s -c, \label{(e5)}\ee
  \be   \dot
p=( \rho-m-\phi \frac{{c}}{s})p, \label{(e6)} \ee and the growth
rate of consumption,with the aid of equations
(\ref{(e4)})-(\ref{(e6)}), is given by  \be \frac{\dot
c}{c}=\phi(\frac{1}{\sigma}-1)m+\phi
\frac{c}{s}+\frac{m-\rho}{\sigma}. \label{(e7)}\ee The initial and
transversality conditions are of the following form:\be c(0)=c_0,
\;s(0)=s_0, \ee and  \be \lim_{t\to\infty} e^{-\rho t} p(t) s(t)=0
\label{(bc)}. \ee
 The partial
Hamiltonian operator determining equation with the aid of equations
(\ref{(e3)})-(\ref{(e6)}) yields
\begin{eqnarray} c^{-\sigma} s^{\phi(1-\sigma)}[\eta_t+(m s -c) \eta_s]-\eta[c^{1-\sigma}\phi s^{\phi(1-\sigma)-1}+c^{-\sigma} s^{\phi(1-\sigma)} m]\nonumber\\
-\bigg[\frac{(cs^{\phi})^{1-\sigma}}{1-\sigma} +c^{-\sigma}
s^{\phi(1-\sigma)}(m s -c)\bigg][\xi_t+(m s -c)
\xi_s]\nonumber\\
=B_t+(m s -c)B_s-\bigg[\eta- \xi(m s-c)\bigg]\rho c^{-\sigma}
s^{\phi(1-\sigma)},\label{e8}\end{eqnarray} in which $\xi(t,s),\;
\eta(t,s)$ and $B(t,s)$. Separating equation (\ref{e8}) with respect
to powers of $c$ yields following overdetermined system for
$\xi(t,s),\; \eta(t,s)$ and $B(t,s)$:
\begin{eqnarray}
c^{2-\sigma}:& \xi_s=0,\nonumber \\
c^{1-\sigma}:&\eta_s+\frac{\phi}{s}\eta+\sigma \xi_t-\rho \xi=0,\nonumber \\
c^{-\sigma}: &\eta_t+m s \eta_s-m \eta - m s \xi_t+\rho(\eta-m s \xi)=0,\label{e9} \\
c:&B_s=0,\nonumber\\
c^0:&B_t=0.\nonumber
\end{eqnarray}
The solution of system (\ref{e9}) provides following partial
Hamiltonian operators $\xi(t,s),\; \eta(t,s)$ and gauge term
$B(t,s)$:
\begin{eqnarray}
\xi=a_1 e^{- \rho t}+a_2 e^{\frac{(\rho-m \phi-m)(1-\sigma)}{\sigma}t},\nonumber\\
\eta=\frac{a_1 s \rho}{(1-\sigma)(\phi+1)}e^{- \rho t}+a_2 m s e^{\frac{(\rho-m \phi-m)(1-\sigma)}{\sigma}t}+a_3 s^{-\phi} e^{-(\rho-m \phi-m)t}, \label{e10}\\
B=a_4. \nonumber
\end{eqnarray}
The first integrals corresponding to operators and gauge terms given
in equation (\ref{e10}) can be computed from formula and are given
as follows:
\begin{eqnarray}
I_1=\frac{ \rho p s e^{-\rho t}}{(\phi+1)(1-\sigma)}-e^{-\rho t}[\frac{(cs^{\phi})^{1-\sigma}}{1-\sigma} +p(m s -c)],\nonumber\\
I_2= e^{\frac{(\rho-m \phi-m)(1-\sigma)}{\sigma}t}[p c-\frac{(cs^{\phi})^{1-\sigma}}{1-\sigma} ], \label{e11}\\
I_3=p s^{-\phi} e^{(m \phi +m-\rho)t}. \nonumber
\end{eqnarray}

\subsubsection{Closed-form solution}
One can utilize either any one or any two of these first integrals
to derive the closed-form solution of system
(\ref{(e4)})-(\ref{(e6)}). Setting $I_3=a$ yields \be
c=a^{-\frac{1}{\sigma}} s^{-\phi}e^{\frac{(m \phi
+m-\rho)}{\sigma}t}.\label{e12}\ee Equation (\ref{e12}) transforms
equation (\ref{(e5)}) into Bernoulli's equation for variable $s$
which gives \be s(t)=\bigg(\frac{\sigma
a^{-\frac{1}{\sigma}}(\phi+1)e^{\frac{m(\phi+1)-\rho}{\sigma}t}}{\rho+m(\sigma-1)(\phi+1)}+b
e^{m(\phi+1)t}\bigg)^{\frac{1}{\phi+1}},\label{e13}\ee and then
equation (\ref{(e5)}) gives $c(t)$.  It is worthy to mention her
that one can directly arrive at expression for $s$ given in
(\ref{e13}) by setting $I_1=b$ and then no integartion is required.
The transversality condition (\ref{(bc)}) takes following form \be
\lim_{t\to\infty} e^{-\rho t} p(t) s(t)=a s^{\phi+1}e^{-m(\phi+1)t}
\label{(bc2)}, \ee and it goes to zero provided
$\rho+m(\sigma-1)(\phi+1)>0$ and $b=0$. The final form of
closed-form solutions subject to initial conditions $c(0)=c_0$ and
$s(0)=s_0$ is given as follows:
\begin{eqnarray}
s(t)=s_0 e^{\frac{m(\phi+1)-\rho}{\sigma(\phi+1)}t},\nonumber\\
c(t)=c_0 e^{\frac{m(\phi+1)-\rho}{\sigma(\phi+1)}t}, \label{e15}\\
p(t)=c_0^{-\sigma}s_0^{\phi(1-\sigma)}e^{(\rho-m-\phi
\frac{c_0}{s_0})t}, \nonumber
\end{eqnarray}
where $\frac{c_0}{s_0}\sigma(\phi+1)=\rho+m(\sigma-1)(\phi+1).$ The
growth rates of all variables of economy are \begin{eqnarray}
\frac{\dot s}{s}= \frac{m(\phi+1)-\rho}{\sigma(\phi+1)},\\
\frac{\dot c}{c}=\frac{m(\phi+1)-\rho}{\sigma(\phi+1)}, \label{e16}\\
\frac{\dot p}{p}=\rho-m-\phi \frac{c_0}{s_0}.
\end{eqnarray}
The variables $s$ and $c$ grow at the same constant rate which makes
true economics sense. The variable $p$ grows at a constant rate
$\rho-m-\phi \frac{c_0}{s_0}$.
\subsection{First integrals for equations of motion for a mechanical system}
Next it is shown how to construct first integrals with the aid of
partial Hamiltonian approach \cite{naz} for more variables case.

A  Hamiltonian function for a mechanical system with non-potential
forces $\Gamma^1=-p_2$ and $\Gamma^2=0$ is
 \be H= \frac{p_1^2}{2}+\frac{p_2^2}{2}+\frac{q_2^2}{2}. \label{rn33}\ee

The equations of motion for a mechanical system are
 \begin{eqnarray} \dot q_1 =p_1,\label{rn34}\\
\dot q_2=p_2,\label{rn35}\\
  \dot p_1
=-p_2,\label{rn36}\\
\dot p_2=-q_2. \label{rn37}\end{eqnarray}  The partial Hamiltonian
operator determining equation {(\ref{(9)})} with the aid of
equations (\ref{rn34})-(\ref{rn37}) yields
\begin{eqnarray} p_1(\eta^1_t+p_1 \eta^1_{q_1}+p_2 \eta^1_{q_2})+ p_2(\eta^2_t+p_1 \eta^2_{q_1}+p_2 \eta^2_{q_2})\nonumber\\
-\eta^2 q_2 -(\frac{p_1^2}{2}+\frac{p_2^2}{2}+\frac{q_2^2}{2})(
p_1\xi_t+p_1 \xi_{q_1}+p_2 \eta^1_{q_2})\label{rn38}\\
=B_t+p_1 B_{q_1}+p_2 B_{q_2}+(\eta^1-\xi
p_1)p_2\nonumber\end{eqnarray} in which $\xi(t,q_1,q_2),\;
\eta^1(t,q_1,q_2),\; \eta^2(t,q_1,q_2)$ and $B(t,q_1,q_2)$.
Separating equation (\ref{rn38}) with respect to powers of $p_1$ and
$p_2$ as $\eta^1,\;\eta^2,\; \xi,\; B$ do not contain $p_1,\; p_2$,
an over determined system for $\eta^1,\;\eta^2,\; \xi,\; B$ is
obtained. After lengthy calculations, following expressions for
$\eta^1,\;\eta^2,\; \xi,\; B$ are obtained:
\begin{eqnarray} \xi=c_1,\; \eta^1=-c_1q_2+c_2 t+c_3, \nonumber\\
\eta^2=c_4 \sin t+c_5 \cos t+c_2,\label{rni}\\
B=\frac{c_1}{2}q_2^2+c_2(q_1-t q_2)-c_3 q_2+c_4 q_2 \cos t -c_5 q_2
\sin t.\nonumber\end{eqnarray} After utilizing formula (\ref{(10)})
the first integrals associated with partial Hamiltonian operators
and gauge terms given in (\ref{rni}) can be expressed as following
form:
\begin{eqnarray} I_1=q_2^2+q_2 p_1+\frac{p_1^2}{2}+\frac{p_2^2}{2},\\
I_2=q_1-t q_2-t p_1-p_2,\\
I_3=q_2+p_1,\\
I_4=q_2 \cos t -p_2 \sin t,\\
I_5= q_2 \sin t +p_2 \cos t.\end{eqnarray} One can utilize these
first integrals to construct closed-form solution for the partial
Hamiltonian system (\ref{rn34})-(\ref{rn37}).
 \subsection{Force-free Duffing Van Der Pol Oscillator}
The standard form of force-free Duffing Van der pol oscillator
equation is
 \be \ddot q+(\alpha+\beta q^2)\dot q-\gamma q+q^3=0, \label{rn1}\ee
 where $\alpha,\; \beta$ and $\gamma$ are arbitrary parameters.
 Equation  (\ref{rn1}) arises in a model describing the propagation
 of voltage pulses along a neuronal axon \cite{osc1,osc2}. Chandrasekar et al \cite{osc1} derived only one
 first integral for the case when $\alpha=\frac{4}{\beta}$
 and $\gamma=-\frac{3}{\beta^2}$.  After making a series of variable transformations and applying the Preller-
Singer method, one first integral was computed in \cite{osc2} and
then inverse transformations provide the first integrals of the
original equations satisfying parameter restriction  $\beta^2
\gamma+3 \alpha \beta -9=0$. I derive here more first integrals of
force-free Duffing Van der pol oscillator using the partial
Hamiltonian approach and these are not reported before in
literature.
 The partial Hamiltonian function for second-order ODE
(\ref{rn1}) is \be H=
-\frac{p^2}{2}+\frac{\gamma}{2}q^2-\frac{q^4}{4} \label{rn2}\ee and
it satisfies

\be \dot q=-p,\label{rn3}\ee \be \dot p=-\gamma q+ q^3+\Gamma
\label{rn4}\ee where $\Gamma=-(\alpha+\beta q^2)p$. The partial
Hamiltonian operator determining equation {(\ref{(9)})} with the aid
of equations (\ref{rn2})-(\ref{rn4})  yields
\begin{eqnarray} p(\eta_t-p \eta_q)-\eta(\gamma
q-q^3)-(-\frac{p^2}{2}+\frac{\gamma}{2}q^2-\frac{q^4}{4})(\xi_t-p
\xi_q)\nonumber\\
=B_t-p B_q+(\eta+p \xi)(\alpha+\beta q^2)p,\label{rn5}\end{eqnarray}
in which $\xi(t,q),\; \eta(t,q)$ and $B(t,q)$. Separating equation
(\ref{rn5}) with respect to powers of $p$ as $\eta, \xi, B$ do not
contain $p$ yields
\begin{eqnarray}
p^3:& \xi_q=0,\nonumber \\
p^2:&-\eta_q+\frac{1}{2}\xi_t=\xi(\alpha+\beta q^2),\nonumber \\
p: &\eta_t+\xi_q(\frac{\gamma}{2}q^2-\frac{q^4}{4})=-B_q+\eta(\alpha+\beta q^2),\label{rn6} \\
p^0:& B_t+\eta(\gamma
q-q^3)+\xi_t(\frac{\gamma}{2}q^2-\frac{q^4}{4})=0.\nonumber
\end{eqnarray}
After lengthy calculations finally system (\ref{rn6}) gives
\begin{eqnarray}
\xi=c_1 e^{\frac{6}{\beta}t},\;\eta=-\frac{1}{\beta}q(\beta^2
q^2+3\alpha \beta-9)c_1e^{\frac{6}{\beta}t}+c_2e^{\frac{3}{\beta}t},\nonumber\\
B=-\frac{q^2}{18 \beta^2}\bigg(\beta^4q^4+(6\alpha
\beta^3-\frac{45}{2}\beta^2)q^2+9\alpha^2 \beta^2-81 \alpha
\beta+162\bigg)c_1 e^{\frac{6}{\beta}t}\label{rn7}\\
+\frac{q}{3\beta}(\beta^2q^2+3\alpha
\beta^2-9)c_2e^{\frac{3}{\beta}t}\nonumber
\end{eqnarray}
provided $\beta^2 \gamma+3 \alpha \beta -9=0$. The first integrals
can be constructed from formula (\ref{(10)}).
 The partial Hamiltonian operators and first integrals are given by
\begin{eqnarray}
X_1=e^{\frac{6}{\beta}t}\frac{\partial }{\partial
t}-\frac{q}{\beta}(\beta^2 q^2+3\alpha \beta-9)e^{\frac{6}{\beta}t}
\frac{\partial }{\partial
q},\nonumber\\
I_1=\frac{1}{2}[p-(\alpha
\beta-3)\frac{q}{\beta}-\frac{\beta}{3}q^3]^2e^{\frac{6}{\beta}t},\nonumber\\
 X_2= e^{\frac{3}{\beta}t}\frac{\partial }{\partial
q},\nonumber\\
I_2=[p-(\alpha
\beta-3)\frac{q}{\beta}-\frac{\beta}{3}q^3]e^{\frac{3}{\beta}t}.\label{rn9}
\end{eqnarray}
 provided $\beta^2 \gamma+3 \alpha \beta -9=0$. Notice that
 $I_1=\frac{I_2^2}{2}$. Now I utilize $I_2$ to construct solution
 of (\ref{rn3}) and (\ref{rn4}). Set $I_2=a_1$ gives
 \be p=(\alpha \beta-3)\frac{q}{\beta}+\frac{\beta q^3}{3}+a_1 e^{-\frac{3}{\beta}t}. \label{ri10}\ee
 Equation (\ref{rn3}) with the aid of equation (\ref{ri10}) results
 in
 \be \dot q(t)+(\alpha \beta-3)\frac{q}{\beta}+\frac{\beta q^3}{3}+a_1 e^{-\frac{3}{\beta}t}=0,\label{rn10} \ee
 where $a_1$ is an arbitrary constant and one can choose it as zero without loss of generality.
  Equation (\ref{rn10}) reduces to a Bernoulli's equation if
  $a_1=0$ which yields
  \be q(t)=\pm \frac{ \sqrt{9a_2(\alpha \beta-3)^2e^{-\frac{2(\alpha \beta-3)}{\beta}t}
  -3 \beta^2(\alpha \beta-3)e^{-\frac{4(\alpha \beta-3)}{\beta}t} }}{\beta^2 e^{-\frac{2(\alpha \beta-3)}{\beta}t} -3 a_2 (\alpha \beta-3)}.\label{ri11}\ee
The exact solution for second order ODE (\ref{rn1}) or its
equivalent system  (\ref{rn3})-(\ref{rn4}) is given in (\ref{ri10})
and (\ref{ri11}). This solution holds provided $\beta^2 \gamma+3
\alpha \beta -9=0$. This solution is new in literature and partial
Hamiltonian approach has made it possible to derive this solution.

\subsection{Lotka-Volterra system}
The notion of partial Hamiltonian function for a dynamical system of
two first-order ODEs is explained with the help of Lotka-Volterra
system \cite{lotka,vol,lot}. Moreover, the first integrals and
closed-form solutions are derived by newly developed partial
Hamiltonian approach \cite{naz}. The two species Lotka-Volterra
model for predator-prey interaction is governed by following two
first-order ODEs:
 \begin{eqnarray}
 \dot q=aq-b pq,\label{rn11}\\
 \dot p= -m p+n p q,\label{rn1111}
 \end{eqnarray}
 where $q$ is the number of prey, $p$ is the number of predator,
 $a$ and $m$ are their per-capita rates
of change in the absence of each other and $b$ and $n$ their
respective rates of change due to interaction.
 A partial Hamiltonian function for system (\ref{rn11})-(\ref{rn1111})
 \be H= a pq-\frac{b}{2}p^2 q \label{rn12}\ee
satisfies \be  \dot q =\frac{\partial H}{\partial p},\label{rn13}\ee
 \be  \dot p =-\frac{\partial H}{\partial q}+\Gamma,\ee
with function $\Gamma$ \be \Gamma=-\frac{b
p^2}{2}+(a-m)p+npq.\label{rn14} \ee The partial Hamiltonian
operators determining equation {(\ref{(9)})} with the aid of
equations (\ref{rn12})-(\ref{rn14}) yields
\begin{eqnarray} p[\eta_t+(aq-bpq) \eta_q]-\eta(ap-\frac{b p^2}{2})-(apq-\frac{b p^2
q}{2})[\xi_t+(a q-b pq)
\xi_q]\nonumber\\
=B_t+(a q-b pq) B_q-[\eta- \xi(a q-b pq)][-\frac{b
p^2}{2}+(a-m)p+npq],\label{rn15}\end{eqnarray} in which $\xi(t,q),\;
\eta(t,q)$ and $B(t,q)$. Separating equation (\ref{rn15}) with
respect to powers of $q$ as $\eta, \xi, B$ do not contain $q$,
provides
\begin{eqnarray}
p^3:& q \xi_q+\xi=0,\nonumber \\
p^2:&\eta_q-\frac{1}{2}\xi_t-\frac{3aq}{2}\xi_q-(\frac{3a}{2}-m+n q)\xi=0,\label{rn16} \\
p: &\eta_t+aq \eta_q+(n q-m)\eta-a q\xi_t-a^2 q^2\xi_q\nonumber \\
&+b qB_q-(a^2 q-a m q+a q^2 n)\xi=0, \nonumber \\
p^0:& B_t+qB_q=0.\nonumber
\end{eqnarray}
The solution of system (\ref{rn16}) yields
\begin{eqnarray}
\xi=\frac{F_1(t)}{q},\nonumber \\
\eta=\bigg(\frac{1}{2} F_1'(t)-m F_1(t)\bigg)\ln(q)+n qF_1(t)+F_2(t),\nonumber \\
B=\bigg(-\frac{1}{4 b}F_1''(t)+\frac{3 m}{4b} F_1'(t)-\frac{m^2}{2b}F_1(t)\bigg)(\ln(q))^2\nonumber\\
+\frac{1}{b}\bigg(-F_2'(t)+m F_2(t)+\frac{(a- n q)}{2}F_1'(t)+m n q F_1(t)\bigg)\ln(q)\nonumber\\
-\frac{n^2 q^2}{2
b}F_1(t)-\frac{n}{b}qF_2(t)-\frac{n}{2b}F'_1(t)q+F_3(t),\label{rn17}
\nonumber
\end{eqnarray}
with $F_1(t)$, $F_2(t)$ and $F_3(t)$ satisfying
\begin{eqnarray}
\frac{1}{4} F_1(t)'''-\frac{3 m}{4}F_1''(t)+\frac{m^2}{2}F_1'(t)=0, \nonumber\\
F_2''(t)-m F_2'(t)-\frac{3 a m}{2}F_1'(t)+a m^2F_1(t)=0, \nonumber\\
\frac{1}{2}F_1''(t)-(m  -\frac{ a}{2}) F_1'(t)- a  m
F_1(t)=0 ,\label{rn18} \\
F_1'(t)+ 2 a F_1(t)=0,\nonumber\\
F_2'(t)+2 a F_2(t)+F_1''(t)+2 a F_1'(t)-2 am
F_1(t)=0,\nonumber\\
F_3'(t)+\frac{a}{b}\bigg(\frac{a}{2}F_1'(t)-F_2'(t)+ m
F_2(t)\bigg)=0.\nonumber
\end{eqnarray}
System (\ref{rn18}) yields \be F_1(t)=C_1 e^{-2 a t},\; F_2(t)=a C_1
e^{-2 a t}-\frac{1}{a}C_2e^{-a t},\; F_3(t)=C_3,\label{rn19}\ee
provided $m=-a$. The partial Hamiltonian operators, gauge term and
first integrals are given by
\begin{eqnarray}
\xi=\frac{e^{-2at}}{q},\; \eta=(a+n q)e^{-2at}, B=-\frac{n^2 q^2}{2
b}e^{-2at},\nonumber\\
I_1=\frac{1}{2b}e^{-2at}(b p+n q)^2,\label{rn19}\\
\xi=0,\; \eta=-\frac{1}{a}e^{-at}, B=-\frac{n q}{a2
b}e^{-at},\nonumber\\
I_2=-\frac{1}{a b}e^{-at}(b p+n q),\nonumber
\end{eqnarray}
provided $a=-m$. Notice that $I_1=\frac{1}{2}a^2 b I_2^2$ and thus
the first integrals are functionally dependent. One can construct
solution of system (\ref{rn11})-(\ref{rn1111}) with the help of one
of these first integrals.

Setting $I_2=\alpha_1$ results in\be -\frac{1}{a b}e^{-at}(b p+n
q)=\alpha_1 \ee
 and this yields
 \be p=-(a \alpha_1 e^{a t} +\frac{n q}{b}) \label{rn20}, \ee
 where $\alpha_1$ is arbitrary constant. Using value of $p$ from equation
 (\ref{rn20}) in equation (\ref{rn11}) yields
 \be \dot q-a(1+b c_1 e^{a t}) q=n q^2,\ee
 and this provides
 \be q(t)=-\frac{a b \alpha_1 e^{a t}}{n-a b \alpha_1 \alpha_2 e^{-b \alpha_1 e^{a t}}}.\label{rn21}\ee
 Equation (\ref{rn20}) and (\ref{rn21}) gives following final form
 of exact solution for variable $p(t)$:
 \be p=-a \alpha_1 e^{a t} +\frac{a n \alpha_1 e^{a t}}{n-a b \alpha_1 \alpha_2 e^{-b \alpha_1 e^{a t}}}. \label{rn22} \ee
 The exact solutions for $q(t)$ and $p(t)$ given in
 equations(\ref{rn21})-(\ref{rn22}) are valid for $a=-m$ case only. The assumption $a=-m$ shows that the per-capita
 rates of change are same for the predator and prey in the absence of each
 other. The solution derived here is not reported before in
 literature and is established due to partial Hamiltonian approach.

\section{Conclusions}
The partial Hamiltonian systems arise widely in different fields of
the applied mathematics. The partial Hamiltonian systems appear for
the mechanical systems with non-holonomic nonlinear constraints and
non-potential generalized forces. In dynamic optimization problems
of economic growth theory involving a non-zero discount factor the
partial Hamiltonian systems arise and are known as the current value
Hamiltonian systems. These systems have a real physical structure.
The partial Hamiltonian approach proposed earlier for a current
value Hamiltonian systems arising in economic growth theory is
applicable to mechanics and other areas.

 In order to show effectiveness of approach
the method is applied to four models: optimal growth model with
environmental asset, the equations of motion for a mechanical
system, the force-free duffing Van Der Pol oscillator and
Lotka-Volterra model. The first integrals and closed-from solutions
of the optimal growth model with environmental asset are derived.
This model yields a current value Hamiltonian system.  The partial
Hamiltonian approach provided three first integrals for this model
and then closed-from solutions are derived with the aid of these
first integrals.  The first integrals of the equations of motion for
a mechanical system with non-potential forces are derived. A
standard Hamiltonian exists for this model and I obtained five first
integrals.  The system of second-order ODEs, describing the
force-free duffing Van Der Pol oscillator, is expressed as a partial
Hamiltonian system of first-order ODEs. The partial Hamiltonian
approach provided two first integrals. Only one of these first
integrals is linearly independent and then this first integral is
utilized to derive closed-form solution. Finally, I derived the
first-integrals and closed-form solutions for Lotka-Volterra system
described by two first-order ODEs. According to best of my knowledge
all these solutions are new in literature and are obtained with the
help of partial Hamiltonian approach.

 It is worthy to mention here that the
partial Hamiltonian systems arise in different fields of applied
mathematics. In each field the functions $\Gamma^i(t,q^i,p_i)$ are
interpreted in different ways but are always some physical
quantities.  The partial Hamiltonian systems arise widely in
non-linear mechanics. In mechanics the functions $\Gamma^i$ contains
non-potential generalized forces or generalized constrained forces
or both. The functions $\Gamma^i$ are associated with discount
factor for economic growth theory.  In other fields as well the
functions $\Gamma^i(t,q^i,p_i)$ have structural properties.
 

\begin{thebibliography}{00}
 \bibitem{naz} Naz, R., Mahomed, F. M., \& Chaudhry, A. (2014). A partial Hamiltonian approach for
current value Hamiltonian systems. Communications in Nonlinear
Science and Numerical Simulation, 19(10), 3600-3610.
 \bibitem{lag}de Lagrange, J. L. (1788). M\'echanique analitique. Paris: Desaint, 1788; 512 p.; in 8.; DCC. 4.403, 1.
 \bibitem{ham}Hamilton, W. R. (2000). The Mathematical Papers of Sir William Rowan Hamilton. CUP Archive.
 \bibitem{leg}Legendre, A. M. (1833).R{\'e}flexions sur differ{\'e}ntes mani{\`e}res de d{\'e}mostrer la th{\'e}orie des parall{\'e}les ou le th{\'e}or{\`e}me sur la somme des
 trois angles du triangle. M{\'e}moires de l'Academie des Sciences de Paris, 13, 213-220.
\bibitem{leg2} Arnold VI (1989) Mathematical methods of classical mechanics.
Springer-Verlag, New York M{\'e}moires de l'Academie des Sciences de
Paris, 13, 213-220.
 \bibitem{eul}Euler, L. (1750). Discovery of a
new principle of mechanics. M{\'e}moires de l'Academie Royale des
science.
\bibitem{bel}Bellman, R. E. (1953). An introduction to the theory of dynamic
programming.
\bibitem{pont}Pontryagin, L. S. (1987). Mathematical theory of optimal processes. CRC Press.
 \bibitem{nlm1}Kgatle, M. R. R.,\& Mason, D. P. (2016). Linear hydraulic fracture with tortuosity: Conservation laws and fluid extraction. International Journal of Non-Linear Mechanics, 79, 10-25.
   \bibitem{nlm2} Kovacic, I. (2006). Conservation laws of two coupled non-linear oscillators, International Journal of Non-Linear Mechanics 41,751-760.
   \bibitem{nlm3}Mahomed, K. S., \& Momoniat, E. (2014). Symmetry classification of first integrals for scalar dynamical equations. International Journal of Non-Linear Mechanics, 59, 52-59.
 \bibitem{Noe}  E. Noether,  Invariante Variationsprobleme.
 Nachr. K\"onig. Gesell. Wissen., G\"ottingen, Math.-Phys. Kl. 1918, Heft 2: 235-257.
(English translation in  Transport Theory and Statistical Physics
1(3) 1971 186-207.)
\bibitem{imran}  Kara, A. H., Mahomed, F. M., Naeem, I., \& Wafo Soh, C. (2007). Partial Noether operators
and first integrals via partial Lagrangians. Mathematical Methods in
the Applied Sciences, 30(16), 2079-2089.
 \bibitem{kar3} Kara, A. H., \& Mahomed, F. M. (2006). Noether-type symmetries and conservation laws via partial Lagrangians. Nonlinear Dynamics, 45(3-4), 367-383.
 \bibitem{cl}Naz, R., Mahomed, F. M., \& Chaudhry, A (2016). A partial Lagrangian Method for Dynamical Systems,  Nonlinear
 dynamics,
 84, 1783–1794.
  \bibitem{nae} Naeem, I., \& Mahomed, F. M. (2009). Approximate partial Noether
operators and first integrals for coupled nonlinear oscillators.
Nonlinear Dynamics, 57(1-2), 303-311.
 \bibitem{ste}H. Steudel (1962). {\"U}ber die Zuordnung zwischen lnvarianzeigenschaften und Erhaltungss{\"a}tzen, Zeitschrift f{\"u}r Naturforschung A, 17(2), 129-132.
 \bibitem{olv} Olver, P. J. (1993). Applications of Lie groups to differential equations, Springer-Verlag, New York.
 \bibitem{anc}Anco, S. C.,  Bluman, G. (2002). Direct construction method for conservation laws of partial differential equations Part I:
 Examples of conservation law classifications. European Journal of Applied Mathematics, 13(05), 545-566.
  \bibitem{gem1}Cheviakov, A. F. (2007). GeM software package for computation of symmetries and conservation laws of differential equations. Computer physics communications, 176(1), 48-61.
  \bibitem{gem2}Cheviakov, A. F. (2010). Computation of fluxes of conservation laws. Journal of Engineering Mathematics, 66(1-3), 153-173.
  \bibitem{comp1}Naz, R., Mahomed, F. M., \& Mason, D. P. (2008). Comparison of different approaches to conservation laws for some partial
differential equations in fluid mechanics. Applied Mathematics and
Computation, 205(1), 212-230.
\bibitem{comp2}Naz, R., Freire, I. L.,
\& Naeem, I. (2014). Comparison of Different Approaches to Construct
First Integrals for Ordinary Differential Equations. In Abstract and
Applied Analysis (Vol. 2014). Hindawi Publishing Corporation.
\bibitem{rom}  Dorodnitsyn, V.,   Kozlov, R., 2010.  Invariance and first
integrals of continuous and discrete Hamitonian equations. Journal
of Engineering Mathematics 66, 253-270.
  \bibitem{naz2016} Naz, R., Chaudhry, A., \& Mahomed, F. M. (2016). Closed-form solutions for the Lucas-Uzawa model of economic growth
  via the partial Hamiltonian approach. Communications in Nonlinear Science and Numerical Simulation, 30(1), 299-306.
   \bibitem{diff}Mahomed, F. M. and Roberts, J. A. G. (2015). Characterization of Hamiltonian symmetries and their first integrals. International Journal of Non-Linear Mechanics, 74, 84-91.
  \bibitem{chiang} Chiang,  A. C. (1992). Elements of dynamic
 optimization. McGraw Hill, New York
  \bibitem{mec}Rongwan, L., \& Jingli, F. (1999). Lie symmetries and conserved quantities of nonconservative nonholonomic systems in phase space. Applied Mathematics and Mechanics, 20(6), 635-640.
    \bibitem{mec1}Fu, J. L., \& Chen, L. Q. (2004). Form invariance, Noether symmetry and Lie symmetry of Hamiltonian systems in phase space. Mechanics Research Communications, 31(1), 9-19.
 \bibitem{erm}Ray, J. R., \& Reid, J. L. (1979). More exact invariants for the time-dependent harmonic oscillator. Physics Letters A, 71(4), 317-318.
  \bibitem{ker}Kermack, W. O., \& McKendrick, A. G. (1927). A contribution to the mathematical theory of epidemics. In Proceedings of the Royal Society of London A: Mathematical,
  Physical and Engineering Sciences (Vol. 115, No. 772, pp. 700-721). The Royal Society.
 Macroeconomic Dynamics, 11(02), 272-289.
 \bibitem{lek}Le Kama, A. A., \& Schubert, K. (2007). A note on the consequences
of an endogenous discounting depending on the environmental quality.
Macroeconomic Dynamics, 11(02), 272-289.
    \bibitem{osc1}Chandrasekar, V. K., Senthilvelan, M., \& Lakshmanan, M. (2005). On the complete integrability and linearization of certain
  second-order nonlinear ordinary differential equations.
 In Proceedings of the Royal Society of London A: Mathematical, Physical and Engineering Sciences (Vol. 461, No. 2060, pp. 2451-2477). The Royal Society.
 \bibitem{osc2}Gao, G., \& Feng, Z. (2010). First integrals for the Duffing-van der Pol type oscillator. Electron. J. Differ. Equations, 19, 1-12.
 \bibitem{lotka}Lotka, A. J. (1925). Elements of physical biology. Williams
\& Wilkins, Baltimore, Maryland, USA.
 \bibitem{vol}Volterra, V. (1931). Variations and Fluctuations of The Number of Individuals in Animal Species Living Together, Translated from 1928 edition by RN Chapman, Animal ecology.
  \bibitem{lot}Berryman, A. A. (1992). The orgins and evolution of predator-prey theory. Ecology, 1530-1535.
  \end{thebibliography}
\end{document}